\documentclass[a4paper]{amsart}
\usepackage{amssymb,amsxtra,mathtools,nicefrac}
\usepackage[T1]{fontenc}
\usepackage{lmodern}
\usepackage[utf8x]{inputenc}

\usepackage[pdftitle={A spectral interpretation for the zeros of the Riemann zeta function},
  pdfauthor={Ralf Meyer},
  pdfsubject={Mathematics; MSC 11M26}]{hyperref}
\usepackage[lite]{amsrefs}
\newcommand*{\MRref}[2]{ \href{http://www.ams.org/mathscinet-getitem?mr=#1}{MR \textbf{#1}}}

\renewcommand{\PrintDOI}[1]{\href{http://dx.doi.org/\detokenize{#1}}{doi: \detokenize{#1}}%
  \IfEmptyBibField{pages}{, (to appear in print)}{}}

\BibSpec{book}{%
    +{}  {\PrintPrimary}                {transition}
    +{,} { \textit}                     {title}
    +{.} { }                            {part}
    +{:} { \textit}                     {subtitle}
    +{,} { \PrintEdition}               {edition}
    +{}  { \PrintEditorsB}              {editor}
    +{,} { \PrintTranslatorsC}          {translator}
    +{,} { \PrintContributions}         {contribution}
    +{,} { }                            {series}
    +{,} { \voltext}                    {volume}
    +{,} { }                            {publisher}
    +{,} { }                            {organization}
    +{,} { }                            {address}
    +{,} { \PrintDateB}                 {date}
    +{,} { }                            {status}
    +{}  { \parenthesize}               {language}
    +{}  { \PrintTranslation}           {translation}
    +{;} { \PrintReprint}               {reprint}
    +{.} { }                            {note}
    +{.} {}                             {transition}
    +{,} { \PrintDOI}                   {doi}
    +{}  {\SentenceSpace \PrintReviews} {review}
}

\usepackage{microtype}

\newcommand*{\nb}{\nobreakdash}

\DeclareUnicodeCharacter{189}{{\nicefrac12}}

\newcommand*{\N}{\mathbb N}
\newcommand*{\Z}{\mathbb Z}
\newcommand*{\Q}{\mathbb Q}
\newcommand*{\R}{\mathbb R}
\newcommand*{\C}{\mathbb C}

\newcommand*{\transpose}[1]{\prescript{\mathrm{t}}{}{#1}}
\newcommand*{\abs}[1]{\lvert#1\rvert}
\newcommand*{\norm}[1]{\lVert#1\rVert}
\newcommand*{\braket}[2]{\langle#1,#2\rangle}
\newcommand*{\ooival}[1]{\mathopen]#1\mathclose[}
\newcommand*{\ocival}[1]{\mathopen]#1\mathclose]}
\newcommand*{\coival}[1]{\mathopen[#1\mathclose[}
\newcommand*{\inv}{^{\times\!}}

\newcommand*{\ima}{\textup i}
\newcommand*{\diff}{\textup d}
\newcommand*{\pv}{\textup{pv}}
\newcommand*{\ID}{\textup{id}}
\newcommand*{\IN}{{\smallint\!}}

\newcommand*{\CCINF}{\mathcal D}
\newcommand*{\Sch}{\mathcal S}
\newcommand*{\Hol}{\mathcal O}

\newcommand*{\Hilm}{\mathcal H}
\newcommand*{\Fourier}{\mathfrak F}
\newcommand*{\Primes}{\mathcal P}

\newcommand*{\hot}{\mathbin{\hat\otimes}}
\newcommand*{\defeq}{\mathrel{\vcentcolon=}}
\newcommand*{\congto}{\xrightarrow\cong}
\newcommand*{\boin}{\prec}
\newcommand*{\prto}{\twoheadrightarrow}
\newcommand*{\into}{\rightarrowtail}
\newcommand*{\addconv}{\dagger}

\DeclareMathOperator{\tr}{tr}
\DeclareMathOperator{\RE}{Re}
\DeclareMathOperator{\End}{End}
\DeclareMathOperator{\Aut}{Aut}
\DeclareMathOperator{\ord}{ord}

\theoremstyle{plain}
\newtheorem{theorem}{Theorem}[section]
\newtheorem{proposition}[theorem]{Proposition}
\newtheorem{lemma}[theorem]{Lemma}
\newtheorem{corollary}[theorem]{Corollary}

\theoremstyle{definition}
\newtheorem{definition}[theorem]{Definition}

\begin{document}

\title[Spectral interpretation]{A spectral interpretation for the zeros of the Riemann zeta function}
\author{Ralf Meyer}
\email{rameyer@uni-math.gwdg.de}
\address{Mathematisches Institut, Georg-August Universität Göttingen, Bunsenstraße 3--5, 37073 Göttingen, Germany}


\begin{abstract}
  Based on work of Alain Connes, I have constructed a spectral interpretation
  for zeros of $L$-functions.  Here we specialise this construction to the
  Riemann $\zeta$-function.  We construct an operator on a nuclear Fréchet
  space whose spectrum is the set of non-trivial zeros of~$\zeta$.  We exhibit
  the explicit formula for the zeros of the Riemann $\zeta$-function as a
  character formula.
\end{abstract}



\thanks{Version of \today.  I thank Folkert Tangermann for carefully reading this
  note and pointing out typos.  I corrected these typos and added two clarifying
  footnotes.}

\maketitle

\section{Introduction}
\label{sec:intro}

The purpose of this note is to explain what the spectral interpretation for
zeros of $L$\nb-functions in~\cite{Meyer:Primes_Rep} amounts to in the
simple special case of the Riemann $\zeta$\nb-function.  The
article~\cite{Meyer:Primes_Rep} is inspired by the work of Alain Connes
in~\cite{Connes:Trace_Formula}.  We will construct a nuclear Fréchet
space~$\Hilm^0_-$ and an operator~$D_-$ on~$\Hilm^0_-$ whose spectrum is equal
to the set of non-trivial zeros of the Riemann $\zeta$\nb-function
$\zeta(s)=\sum_{n=1}^\infty n^{-s}$.  By definition, the non-trivial zeros
of~$\zeta$ are the zeros of the \emph{complete $\zeta$\nb-function}
\begin{equation}  \label{eq:complete_zeta}
  \xi(s) = \pi^{-s/2} \Gamma(s/2)\zeta(s).
\end{equation}
In addition, the algebraic multiplicity of~$s$ as an eigenvalue of~$D_-$ is
the zero order of~$\xi$ at~$s$.  Thus~$D_-$ is a spectral interpretation for
the zeros of~$\xi$.

We construct~$D_-$ as the generator of a smooth representation~$\rho_-$ of
$\R\inv_+\cong\R$ on~$\Hilm^0_-$.  Although the single operator~$D_-$ is more
concrete, it is usually better to argue with the representation~$\rho_-$
instead.  Let $\CCINF(\R\inv_+)$ be the convolution algebra of smooth,
compactly supported functions on~$\R\inv_+$.  The integrated form of~$\rho_-$
is a bounded algebra homomorphism $\IN\rho_-\colon
\CCINF(\R\inv_+)\to\End(\Hilm^0_-)$.  We show that $\rho_-$ is a summable
representation in the notation of~\cite{Meyer:Primes_Rep}.  That is,
$\IN\rho_-(f)$ is a nuclear operator for all $f\in\CCINF(\R\inv_+)$.  The
character $\chi(\rho_-)$ is the distribution on~$\R\inv_+$ defined by
$\chi(\rho_-)(f)=\tr \IN\rho_-(f)$.  The representation~$\rho_-$ is part of a
virtual representation $\rho= \rho_+ \ominus \rho_-$, where~$\rho_+$ is a
spectral interpretation for the poles of~$\xi$.  That is, $\rho_+$ is
$2$\nb-dimensional and its generator~$D_+$ has eigenvalues $0$ and~$1$.  We
interpret~$\rho$ as a formal difference of~$\rho_+$ and~$\rho_-$ and therefore
define $\chi(\rho) \defeq \chi(\rho_+)-\chi(\rho_-)$.

The spectrum of~$\rho$ consists exactly of the poles and zeros of~$\xi$, and
the spectral multiplicity (with appropriate signs) of $s\in\C$ is the order
of~$\xi$ at~$s$, which is positive at the two poles $0$ and~$1$ and negative
at the zeros of~$\xi$.  The spectral computation of the character yields
\[
\chi(\rho)(f) = \sum_{s\in\C} \ord_\xi(s) \hat{f}(s),
\qquad
\text{where}
\quad
\hat{f}(s) \defeq \int_0^\infty f(x)x^s \,\frac{\diff x}{x}.
\]

If an operator has a sufficiently nice integral kernel, then we may also
compute its trace by integrating its kernel along the diagonal.  This recipe
applies to $\IN\rho_-(f)$ and yields another formula for $\chi(\rho_-)$.
Namely,
\[
\chi(\rho) = W = \sum_{p\in\Primes} W_p + W_\infty,
\]
where~$\Primes$ is the set of primes,
\begin{align}  \label{eq:Wp}
  W_p(f) &= \sum_{e=1}^\infty f(p^{-e}) p^{-e} \ln(p)
  + \sum_{e=1}^\infty f(p^e) \ln(p),
  \\  \label{eq:Winfty}
  W_\infty(f) &= \pv \int_0^\infty
  \frac{f(x)}{\abs{1-x}} + \frac{f(x)}{1+x}
  \,\diff x.
\end{align}
The distribution~$W_\infty$ involves a principal value because the integrand
may have a pole at~$1$.  Equating the two formulas for $\chi(\rho)$, we get
the well-known explicit formula that relates zeros of~$\xi$ and prime numbers.

We do not need the functions $\zeta$ and~$\xi$ to define our spectral
interpretation.  Instead we use an operator~$Z$, called the Zeta operator,
which is closely related to the $\zeta$\nb-function.  This operator is the
key ingredient in our construction.  In addition, we have to choose the domain
and target space of~$Z$ rather carefully.  It is possible to prove the Prime
Number Theorem using the representation~$\rho$ instead of the
$\zeta$\nb-function.  The only input that we need is that the
distribution~$W$ is quantised, that is, of the form $\sum n(s)\hat{f}(s)$ with
some function $n\colon \C\to\Z$ (see \cite{Meyer:Primes_Rep}).

Our constructions for the $\zeta$\nb-function generalise to
Dirichlet $L$\nb-functions.  We indicate how this is done in the last
section.  We run into problems, however, for number fields with more than one
infinite place.  There are also other conceptual and aesthetic reasons for
prefering adèlic constructions as in \cite{Meyer:Primes_Rep}.  Our goal here
is only to make these constructions more explicit in a simple special case.

\section{The ingredients: some function spaces and operators}
\label{sec:ingredients}

Let $\Sch(\R)$ be the Schwartz space of~$\R$.  Thus $f\colon \R\to\C$ belongs
to $\Sch(\R)$ if and only if all its derivatives~$f^{(n)}$ are rapidly
decreasing in the sense that $f^{(n)}(x) = O(\abs{x}^{-s})$ for
$\abs{x}\to\infty$ for all $s\in\R_+$, $n\in\N$.  We topologise $\Sch(\R)$ in
the usual fashion.  The convolution turns $\Sch(\R)$ into a Fréchet algebra.

We remark that we neither gain nor loose anything if we view $\Sch(\R)$ as a
bornological vector space as in~\cite{Meyer:Primes_Rep}.  All function spaces
that we shall need are Fréchet spaces, so that bornological and topological
analysis are equivalent.  The bornological point of view only becomes superior
if we mix $\Sch(\R)$ with spaces like $\Sch(\Q_p)$, which are not Fréchet.

We use the natural logarithm $\ln$ to identify the multiplicative
group~$\R\inv_+$ with~$\R$.  This induces an isomorphism between the Schwartz
algebras $\Sch(\R\inv_+)$ and $\Sch(\R)$.  The standard Lebesgue measure
on~$\R$ corresponds to the Haar measure $\diff\inv x= x^{-1}\,\diff x$ on~$\R\inv_+$.
We always use this measure in the following.

We let
\[
\Sch(\R\inv_+)_s
= \Sch(\R\inv_+)\cdot x^{-s}
= \{f\colon \R\inv_+\to\C \mid (x\mapsto f(x) x^s) \in\Sch(\R\inv_+)\}
\]
for $s\in\R$ and
\[
\Sch(\R\inv_+)_I = \bigcap_{s\in I} \Sch(\R\inv_+)_s
\]
for an interval $I\subseteq\R$.  We will frequently use that
\[
\Sch(\R\inv_+)_{[a,b]} = \Sch(\R\inv_+)_a\cap \Sch(\R\inv_+)_b.
\]
The reason for this is that $\Sch(\R\inv_+)$ is closed under multiplication
by $(x^\epsilon+x^{-\delta})^{-1}$ for $\epsilon,\delta\ge0$.

Hence $\Sch(\R\inv_+)_{[a,b]}$ becomes a Fréchet space in a canonical way.
Exhausting~$I$ by an increasing sequence of compact intervals, we may turn
$\Sch(\R\inv_+)_I$ into a Fréchet space for general~$I$.  Since $x\mapsto x^s$
is a character of~$\R\inv_+$, the spaces $\Sch(\R\inv_+)_s$ for $s\in\R$ are
closed under convolution.  Hence $\Sch(\R\inv_+)_I$ is closed under
convolution as well and becomes a Fréchet algebra.  We are particularly
interested in
\begin{align*}
  \Hilm_- = \Hol(\R\inv_+) &\defeq \Sch(\R\inv_+)_{\ooival{-\infty,\infty}},
  \\
  \Sch_> &\defeq \Sch(\R\inv_+)_{\ooival{1,\infty}},
  \\
  \Sch_< &\defeq \Sch(\R\inv_+)_{\ooival{-\infty,0}}.
\end{align*}
We also let
\[
\Hilm_+ \defeq \{ f\in\Sch(\R)\mid \text{$f(x)=f(-x)$ for all $x\in\R$}\}.
\]
The spaces $\Hilm_\pm$ will be crucial for our spectral interpretation; the
spaces $\Sch_>$ and $\Sch_<$ only play an auxiliary role as sufficiently large
spaces in which the others may be embedded.

Given topological vector spaces $A$ and~$B$, we write $A\boin B$ to denote
that~$A$ is contained in~$B$ and that the inclusion is a continuous linear
map.  Clearly,
\[
\Sch(\R\inv_+)_I\boin\Sch(\R\inv_+)_J  \qquad\text{if $I\supseteq J$.}
\]

The group~$\R\inv_+$ acts on $\Sch(\R\inv_+)_I$ and~$\Hilm_+$ by the
\emph{regular representation}
\[
\lambda_t f(x) \defeq f(t^{-1}x)  \qquad\text{for $t,x\in\R\inv_+$.}
\]
Its integrated form is given by the same formula as the convolution:
\[
\IN\lambda(h) f(x) \defeq \int_0^\infty h(t) f(t^{-1}x) \,\diff\inv t.
\]

We denote the projective complete topological tensor product by~$\hot$ (see
\cite{Grothendieck:Produits}).  If $V$ and~$W$ are Fréchet spaces,
so is $V\hot W$.  We want to know $\Sch(\R\inv_+)_I\hot\Sch(\R\inv_+)_J$ for
two intervals $I,J$.  Since both tensor factors are nuclear Fréchet spaces,
this is easy enough to compute.  We find
\begin{multline}  \label{eq:weightedSch_tensor}
  \Sch(\R\inv_+)_I\hot\Sch(\R\inv_+)_J
  \\ \cong \{f\colon (\R\inv_+)^2\to\C \mid
  \text{$f(x,y)\cdot x^sy^t \in \Sch((\R\inv_+)^2)$ for all $s\in I$, $t\in
    J$}\}
\end{multline}
with the canonical topology.  This follows easily from
$\Sch(\R)\hot\Sch(\R)\cong\Sch(\R^2)$ and the compatibility of~$\hot$ with
inverse limits.

We shall need the \emph{Fourier transform}
\begin{equation}  \label{eq:def_Fourier}
  \Fourier\colon \Sch(\R)\to\Sch(\R),  \qquad
  \Fourier f(y) \defeq \int_\R f(x) \exp(2\pi\ima xy) \,\diff x.
\end{equation}
Notice that $\Fourier f$ is even if~$f$ is.  Hence~$\Fourier$ restricts to an
operator on $\Hilm_+\subseteq\Sch(\R)$.  In the following, we usually
restrict~$\Fourier$ to this subspace.  It is well-known that $\Fourier^{-1}
f(y)=\Fourier f(-y)$ for all $f\in\Sch(\R)$, $y\in\R$.  Hence
\begin{equation}  \label{eq:Fourier_square}
  \Fourier^2 = \ID \qquad\text{as operators on $\Hilm_+$.}
\end{equation}
Since~$\Fourier$ is unitary on $L^2(\R,\diff x)$, it is also unitary on the
subspace of even functions, which is isomorphic to $L^2(\R\inv_+,x\,\diff\inv x)$.

We also need the involution
\begin{equation}  \label{eq:def_J}
  J\colon \Hol(\R\inv_+)\to\Hol(\R\inv_+),  \qquad
  Jf(x) \defeq x^{-1} f(x^{-1}).
\end{equation}
We have $J^2=\ID$.  One checks easily that~$J$ extends to a unitary operator
on $L^2(\R\inv_+,x\,\diff\inv x)$ and to an isomorphism of topological vector
spaces
\begin{equation}  \label{eq:J_reflects}
  J\colon \Sch(\R\inv_+)_I \congto \Sch(\R\inv_+)_{1-I}
\end{equation}
for any interval~$I$.  Especially, $J$ is an isomorphism between $\Sch_<$
and~$\Sch_>$.

We have
\begin{equation}  \label{eq:Fourier_variant}
  \IN\lambda(h) \circ \Fourier = \Fourier\circ \IN\lambda(Jh),
  \qquad
  \IN\lambda(h) \circ J = J\circ \IN\lambda(Jh)
\end{equation}
for all $h\in\CCINF(\R\inv_+)$.  Hence the composites $\Fourier J$ and
$J\Fourier=(\Fourier J)^{-1}$, which are unitary operators on
$L^2(\R\inv_+,x\,\diff\inv x)$, commute with the regular representation~$\lambda$.

A \emph{multiplier} of $\Sch(\R\inv_+)_I$ is a continuous linear operator on
$\Sch(\R\inv_+)_I$ that commutes with the regular representation.

\begin{proposition}  \label{pro:Sch_estimates}
  Viewing $\Hilm_+\subseteq L^2(\R\inv_+,x\,\diff\inv x)$, we have
  \begin{equation}    \label{eq:Hilmplus_Fourier}
    \Hilm_+ = \{f\in L^2(\R\inv_+,x\,\diff\inv x)\mid
    \text{$f\in \Sch(\R\inv_+)_{\ooival{0,\infty}}$ and
      $J\Fourier(f)\in\Sch(\R\inv_+)_{\ooival{-\infty,1}}$}\}.
  \end{equation}
  \begin{itemize}
  \item The operator $\Fourier J$ is a multiplier of $\Sch(\R\inv_+)_I$ for
    $I\subseteq\ooival{0,\infty}$.
    
  \item The operator $J \Fourier$ is a multiplier of $\Sch(\R\inv_+)_I$ for
    $I\subseteq\ooival{-\infty,1}$.
    
  \item Hence $\Fourier J$ and $J\Fourier$ are invertible multipliers of
    $\Sch(\R\inv_+)_I$ for $I\subseteq\ooival{0,1}$.

  \end{itemize}
\end{proposition}

\begin{proof}
  We first have to describe $\Sch(\R\inv_+)_I$ more explicitly.  For
  simplicity, we assume the interval~$I$ to be open.  Let $Df(x)\defeq x\cdot
  f'(x)$.  This differential operator is the generator of the
  representation~$\lambda$ of $\R\inv_+\cong\R$.  Let $f\in
  L^2(\R\inv_+,x\,\diff\inv x)$.  Then $f\in\Sch(\R\inv_+)_I$ if and only if
  $D^m(f\cdot x^s) \cdot (\ln x)^k\in L^2(\R\inv_+,\diff\inv x)$ for all
  $m,k\in\N$, $s\in I$.  Using the Leibniz rule, one shows that this is
  equivalent to $D^m(f)\cdot x^s \cdot (\ln x)^k\in L^2(\R\inv_+,\diff\inv x)$ for
  all $m,k\in\N$, $s\in I$.  Since~$I$ is open, we may replace~$x^s$ by
  $x^{s+\epsilon}+x^{s-\epsilon}$ for some $\epsilon>0$.  This dominates $x^s
  (\ln x)^k$ for any $k\in\N$, so that it suffices to require $(D^m f)\cdot
  x^s\in L^2(\R\inv_+,\diff\inv x)$ for all $m\in\N$, $s\in I$.
  
  This description of $\Sch(\R\inv_+)_I$ easily implies $\Hilm_+\subseteq
  \Sch(\R\inv_+)_{\ooival{0,\infty}}$.  Since~$\Fourier$ maps~$\Hilm_+$ to
  itself and~$J$ maps $\Sch(\R\inv_+)_{\ooival{0,\infty}}$ to
  $\Sch(\R\inv_+)_{\ooival{-\infty,1}}$, we get ``$\subseteq$''
  in~\eqref{eq:Hilmplus_Fourier}.  Conversely, $f\in\Hilm_+$ if and only
  if~$f$ and $\Fourier f$ are both $O(x^{-s})$ for $\abs{x}\to\infty$ for all
  $s\in\N$.  This yields ``$\supseteq$'' in~\eqref{eq:Hilmplus_Fourier} and
  finishes the proof of~\eqref{eq:Hilmplus_Fourier}.
  
  In the following computation, we describe
  $\Hol(\R\inv_+)\hot\Sch(\R\inv_+)_I$ as in~\eqref{eq:weightedSch_tensor}.
  Choose any $\psi\in\Hol(\R\inv_+)$ with $\int_0^\infty \psi(x)\,\diff\inv x=1$.
  Then $\sigma f(x,y)\defeq \psi(x) f(xy)$ defines a continuous linear map
  from $\Sch(\R\inv_+)_I$ to $\Hol(\R\inv_+)\hot\Sch(\R\inv_+)_I$
  by~\eqref{eq:weightedSch_tensor}.  This is a section for the convolution map
  \[
  \IN\lambda\colon \Hol(\R\inv_+)\hot\Sch(\R\inv_+)_I\to\Sch(\R\inv_+)_I,
  \qquad
  (\IN\lambda f)(x) = \int_0^\infty f(t,t^{-1}x)\,\diff\inv t.
  \]
  Let $I\subseteq \ooival{0,\infty}$.  Then
  $\Hol(\R\inv_+) = J\Hol(\R\inv_+) \subset \Hilm_+ =
  \Fourier\Hilm_+ \subset \Sch(\R\inv_+)_I$ if
  Hence we get a continuous linear operator
  \[
  \Sch(\R\inv_+)_I
  \overset{\sigma}\longrightarrow \Hol(\R\inv_+)\hot\Sch(\R\inv_+)_I
  \overset{\Fourier J\hot\ID}\longrightarrow
  \Sch(\R\inv_+)_I\hot\Sch(\R\inv_+)_I
  \overset{\IN\lambda}\longrightarrow \Sch(\R\inv_+)_I.
  \]
  The last map exists because $\Sch(\R\inv_+)_I$ is a convolution algebra.
  If we plug $f_0*f_1$ with $f_0,f_1\in\Hol(\R\inv_+)$ into this operator, we
  get $\Fourier J(f_0*f_1)$ because~$\sigma$ is a section for $\IN\lambda$ and
  because of~\eqref{eq:Fourier_variant}.  Since products $f_0*f_1$ are dense
  in $\Hol(\R\inv_+)$, the above operator on $\Sch(\R\inv_+)_I$ extends
  $\Fourier J$ on $\Hol(\R\inv_+)$.
  
  Now~\eqref{eq:J_reflects} implies the continuity of $J\Fourier=J(\Fourier
  J)J$ on $\Sch(\R\inv_+)_I$ for $I\subseteq\ooival{-\infty,1}$.  Hence both
  $J\Fourier$ and $\Fourier J$ are multipliers of $\Sch(\R\inv_+)_I$ for
  $I\subseteq\ooival{0,1}$.  They are inverse to each other on
  $\Sch(\R\inv_+)_I$ because they are inverse to each other on
  $L^2(\R\inv_+,x\,\diff\inv x)$.
\end{proof}

\section{The Zeta operator and the Poisson Summation Formula}
\label{sec:Zeta_operator}

In this section we study the properties of the following operator:

\begin{definition}  \label{def:zeta_operator}
  The \emph{Zeta operator} is defined by
  \[
  Zf(x) \defeq \sum_{n=1}^\infty f(nx) = \sum_{n=1}^\infty \lambda_n^{-1} f(x)
  \qquad\text{for }f\in\Hilm_+,\ x\in\R\inv_+.
  \]
\end{definition}

Let~$\check{\zeta}$ be the distribution $\sum_{n=1}^\infty \delta_n^{-1}\colon
\psi\mapsto \sum_{n=1}^\infty \psi(n^{-1})$.  Then
\[
\check{\zeta}\sphat(s) = \sum_{n=1}^\infty n^{-s} = \zeta(s),
\qquad
Zf= \IN\lambda(\check{\zeta})(f).
\]
Thus the data $Z$, $\zeta$, and~$\check{\zeta}$ are equivalent.

The Euler product expansion of the $\zeta$\nb-function takes the following
form in this picture.  Let~$\Primes$ be the set of prime numbers in~$\N^*$.
We have
\[
Z = \prod_{p\in\Primes} \sum_{e=0}^\infty \lambda_{p^{-e}}
= \prod_{p\in\Primes} (1-\lambda_p^{-1})^{-1}.
\]
Hence we get a candidate for an inverse of~$Z$:
\[
Z^{-1} f(x) = \prod_{p\in\Primes} (1-\lambda_p^{-1}) f(x)
= \sum_{n=1}^\infty \mu(n) f(nx).
\]
Here $\mu(n)$ is the usual Möbius function; it vanishes unless~$n$ is
square-free, and is $(-1)^j$ if~$n$ is a product of~$j$ different prime
numbers.

The following assertion is equivalent to the absolute convergence in the
region $\RE s>1$ of the Euler product defining $\zeta(s)$.

\begin{proposition}  \label{pro:Z_invertible}
  $Z$ and~$Z^{-1}$ are continuous linear operators on~$\Sch_>$ that are
  inverse to each other.
\end{proposition}

\begin{proof}
  Let $Df(x)= xf'(x)$.  We have observed above that $f\in\Sch_>$ if and only
  if $D^m f\in L^2(\R\inv_+,x^{2s}\,\diff\inv x)$ for all $m\in\N$, $s>1$.  We
  check that $Z$ and~$Z^{-1}$ preserve this estimate.  The
  operator~$\lambda_t$ for $t\in\R\inv_+$ has norm $\norm{\lambda_t}^s=t^s$ on
  $L^2(\R\inv_+,x^{2s}\,\diff\inv x)$ for all $s\in\R$.  Hence the same estimates
  that yield the absolute convergence of the Euler product for $\zeta(s)$ also
  show that the products $\prod_{p\in\Primes} (1-\lambda_p^{-1})^{\pm1}$
  converge absolutely with respect to the operator norm on
  $L^2(\R\inv_+,x^{2s}\,\diff\inv x)$ for $s>1$.  Since all factors commute
  with~$D$, this remains true if we replace $L^2(\R\inv_+,x^{2s}\,\diff\inv x)$ by
  the Sobolev space defined by the norm
  \[
  \norm{f}^{m,s}_2 \defeq \int_0^\infty \abs{D^m f(x)}^2 x^{2s} \,\diff\inv x.
  \]
  Thus $Z$ and~$Z^{-1}$ are continuous linear operators on these Sobolev
  spaces for all $m\in\N$ and all $s>1$.  This yields the assertion.
\end{proof}

Next we recall the \emph{Poisson Summation Formula}.  It asserts that
\[
\sum_{n\in\Z} f(xn) = x^{-1} \sum_{n\in\Z} \Fourier f(x^{-1}n)
\]
for all $x\in\R\inv_+$ and all $f\in\Sch(\R)$.  For $f\in\Hilm_+$, this
becomes
\begin{equation}  \label{eq:Poisson}
  \frac{f(0)}{2} + Z f(x) = JZ\Fourier f(x) + \frac{\Fourier f(0)}{2x}.
\end{equation}
Let
\[
\Hilm_\cap \defeq \{f\in\Sch(\R)\mid f\text{ even, }f(0)=\Fourier f(0)=0\}.
\]
This is a closed, $\lambda$\nb-invariant subspace of $\Sch(\R)$.  If
$f\in\Hilm_\cap$, then~\eqref{eq:Poisson} simplifies to $Zf= JZ\Fourier f$.

Let $h_0\colon \R\inv_+\to[0,1]$ be a smooth function with $h_0(t)=1$ for
$t\ll1$ and $h_0(t)=0$ for $t\gg1$.  Let~$\Hilm_\cup$ be the space of
functions on $\R\inv_+$ that is generated by $\Hilm_-=\Hol(\R\inv_+)$ and the
two additional functions $h_0$ and $x^{-1}\cdot h_0$.  The regular
representation extends to~$\Hilm_\cup$.  Writing $\C(x^s)$ for~$\C$ equipped
with the representation by the character~$x^s$, we get an extension of
representations
\begin{equation}  \label{eq:pole_extension}
  \Hilm_- \into \Hilm_\cup \prto \C(x^0)\oplus \C(x^1).
\end{equation}

\begin{theorem}  \label{the:Zeta_estimate}
  The Zeta operator is a continuous linear map $Z\colon \Hilm_+\to\Hilm_\cup$.
  Even more, this map has closed range and is a topological isomorphism onto
  its range.  We have $Zf\in\Hilm_-$ if and only if $f\in\Hilm_\cap$.
\end{theorem}

\begin{proof}
  Proposition~\ref{pro:Sch_estimates} yields $\Hilm_+\boin\Sch_>$.  By
  Proposition~\ref{pro:Z_invertible}, $Z$ is continuous on~$\Sch_>$.  Hence we
  get continuity of $Z\colon \Hilm_+\to\Sch_>$.  Similarly, $JZ\Fourier$ is a
  continuous linear operator $\Hilm_+\to\Sch_<$.  By~\eqref{eq:Poisson}, $Z$
  restricts to a continuous linear map
  \[
  \Hilm_\cap\to\Sch_>\cap\Sch_<=\Hilm_-.
  \]
  Equation~\eqref{eq:Poisson} also implies that~$Zf$ still belongs to~$\Hilm_\cup$
  for arbitrary $f\in\Hilm_+$ and that $Zf\in\Hilm_-$ if and only if
  $f\in\Hilm_\cap$.
  
  It remains to prove that~$Z$ is a topological isomorphism onto its range.
  This implies that the range is closed because all spaces involved are
  complete.  It suffices to prove that the restriction $Z\colon
  \Hilm_\cap\to\Hilm_-$ is an isomorphism onto its range.  Equivalently, a
  sequence $(f_n)$ in~$\Hilm_\cap$ converges if and only if $(Zf_n)$ converges
  in~$\Hilm_-$.  One implication is contained in the continuity of~$Z$.
  Suppose that the sequence $(Zf_n)$ converges in~$\Hilm_-$.  Hence it
  converges in both $\Sch_>$ and~$\Sch_<$.  Equation~\eqref{eq:Poisson} yields
  $Zf_n=JZ\Fourier f_n$.  Using~\eqref{eq:J_reflects}, we get that both
  $(Zf_n)$ and $(Z\Fourier f_n)$ converge in~$\Sch_>$.
  Proposition~\ref{pro:Z_invertible} yields that $(f_n)$ and $(\Fourier f_n)$
  converge in~$\Sch_>$.
  
  Therefore, $(D^mf_n\cdot x^s)$ and $(D^m\Fourier f_n\cdot x^s)$ converge in
  $L^2(\R,\diff x)$ for all $s>½$, where $Df(x)=xf'(x)$.  Hence $x^k(\diff/\diff x)^l f_n$
  and $(\diff/\diff x)^k (x^l f_n)$ converge if $k,l\in\N$ satisfy $k>l+1$.  The first
  condition implies convergence in $\Sch(\R\setminus \ooival{-1,1})$ because
  $x\ge1$ in this region.  The second condition contains convergence of
  $(\diff/\diff x)^k f_n$ in $L^2([-1,1],\diff x)$.  Hence both conditions together imply
  convergence in $\Sch(\R)$ as desired.
\end{proof}

We now discuss the close relationship between the above theorem and the
meromorphic continuation of the $\zeta$\nb-function and the functional
equation (see also \cite{Tate:Thesis}).  Recall that $\hat{f}(s)\defeq
\int_0^\infty f(x) x^s\,\diff\inv x$.  This defines an entire function for
$f\in\Hilm_-$.  We have described~$Z$ as the convolution with the
distribution~$\check{\zeta}$ on~$\R\inv_+$, which satisfies
$\check{\zeta}\sphat(s)=\zeta(s)$.  Therefore,
\begin{equation}  \label{eq:Zspectral}
  (Zf)\sphat(s)=\zeta(s)\cdot \hat{f}(s)
  \qquad
  \text{for all }f\in\Sch_>,\ s\in\C\text{ with }\RE s>1.
\end{equation}
The Poisson Summation Formula
implies $Zf\in \Hilm_-$ for $f\in\Hilm_\cap$, so that $\zeta(s)\hat{f}(s)$
extends to an entire function.  Especially, this holds if $f\in\Hilm_-$
satisfies $\hat{f}(1)=0$.  For such~$f$, the function $\hat{f}(s)$ is an
entire function on~$\C$ as well.  Therefore, $\zeta$ has a meromorphic
continuation to all of~$\C$.  For any $s\neq1$, there is $f\in\Hilm_-$
with $\hat{f}(1)=0$ and $\hat{f}(s)\neq0$.  Therefore, the only possible pole
of~$\zeta$ is at~$1$.

It is easy to see that $(Jf)\sphat(s)=\hat{f}(1-s)$.  Hence~\eqref{eq:Poisson}
implies
\begin{equation}  \label{eq:functional_equation}
  \zeta(1-s)\cdot (J\Fourier f)\sphat(s) = \zeta(s)\cdot \hat{f}(s)
  \qquad\text{for all }s\in\C,\ f\in\Hilm_\cap.
\end{equation}
This equation still holds for
$f\in\Hilm_+$ by $\R\inv_+$\nb-equivariance.

Now we plug in the special function $f(x)=2\exp(-\pi x^2)$, which
satisfies $\Fourier f=f$ and $\hat{f}(s)=\pi^{-s/2}\Gamma(s/2)$.  Thus
$\zeta(s)\hat{f}(s)=\xi(s)$ is the complete $\zeta$\nb-function.
Equation~\eqref{eq:functional_equation} becomes the functional
equation $\xi(1-s)=\xi(s)$.

\section{The spectral interpretation}
\label{sec:spectral_interpretation}

Let $Z\Hilm_+\subseteq\Hilm_\cup$ be the range of~$Z$.  This is a closed
subspace of~$\Hilm_\cup$ and topologically isomorphic to~$\Hilm_+$ by
Theorem~\ref{the:Zeta_estimate}.  Moreover,
\[
Z\Hilm_\cap=Z\Hilm_+\cap\Hilm_-,
\qquad
\Hilm_\cup=Z\Hilm_+ + \Hilm_-.
\]
We define
\begin{align*}
  \Hilm^0_+ &\defeq \Hilm_+/\Hilm_\cap \cong \Hilm_\cup/\Hilm_-,
  \\
  \Hilm^0_- &\defeq \Hilm_-/Z\Hilm_\cap \cong \Hilm_\cup/Z\Hilm_+.
\end{align*}
We equip $\Hilm^0_\pm$ with the quotient topology from~$\Hilm_\pm$ or
from~$\Hilm_\cup$ (both topologies on $\Hilm^0_\pm$ coincide) and with the
representations~$\rho_\pm$ of~$\R\inv_+$ induced by~$\lambda$ on $\Hilm_\pm$
or~$\Hilm_\cup$.  We view the pair of representations $(\rho_+,\rho_-)$ as a
formal difference $\rho_+\ominus\rho_-$, that is, as a virtual representation
of~$\R\inv_+$.

A smooth representation of $\R\inv_+\cong\R$ is determined uniquely by the
action of the generator of the Lie algebra of $\R\inv_+$.  This generator
corresponds to the scaling-invariant vector field $Df(x)=xf'(x)$ on
$\Hilm_\cup$.  We let $D_\pm$ be the operators on $\Hilm^0_\pm$ induced by~$D$
on~$\Hilm_\cup$.\footnote{These act by \(f\mapsto -xf'(x)\) because of the
  inverses in the definition of the regular representation~\(\lambda\).}
We claim that the operators~$D_\pm$ are spectral
interpretations for the poles and zeros of the complete
$\zeta$\nb-function~$\xi$ defined in~\eqref{eq:complete_zeta}.  This is
rather trivial for~$D_+$.  Since
\[
\Hilm^0_+ \cong \C(x^0)\oplus\C(x^1),
\]
the operator~$D_+$ is equivalent to the diagonal $2\times2$-matrix
with eigenvalues $0$, $1$.  Hence it is a spectral interpretation for the two
poles of~$\xi$.  To treat~$\Hilm^0_-$, we identify the Fourier--Laplace
transforms of $\Hilm_-$ and $Z\Hilm_\cap$.  Recall that
$\hat{f}(s)=\int_0^\infty f(x) x^s\,\diff\inv x$.

\begin{theorem}  \label{the:Lap_range}
  The operator $f\mapsto \hat{f}$ identifies~$\Hilm_-$ with the space of
  entire functions $h\colon \C\to\C$ for which $t\mapsto h(s+\ima t)$ is a
  Schwartz function on~$\R$ for each $s\in\R$.
  
  The subspace $Z\Hilm_\cap$ is mapped to the space of those entire
  functions~$h$ for which
  \[
  t\mapsto \frac{h(s+\ima t)}{\zeta(s+\ima t)}
  \qquad\text{and}\qquad
  t\mapsto \frac{h(s+\ima t)}{\zeta(1-s-\ima t)}
  \]
  are Schwartz functions for $s\ge½$ and $s\le½$,
  respectively.  \textup{(}In particular, this means that
  $h(z)/\zeta(z)$ has no poles with $\RE z\ge½$ and
  $h(z)/\zeta(1-z)$ has no poles with $\RE z\le½$.\textup{)}
\end{theorem}

\begin{proof}
  It is well-known that the Fourier transform is an isomorphism (of
  topological vector spaces) $\Sch(\R\inv_+)\cong\Sch(\ima\R)$.  Hence
  $f\mapsto \hat{f}$ is an isomorphism $\Sch(\R\inv_+)_s\cong\Sch(s+\ima\R)$
  for all $s\in\R$.  It is clear that $\hat{f}$ is an entire function on~$\C$
  for $f\in\Hilm_-$.  Since $\Hilm_-=\bigcap_{s\in\R} \Sch(\R\inv_+)_s$, we
  also get $\hat{f}\in\Sch(s+\ima\R)$ for all $s\in\R$.  Conversely, if~$h$ is
  entire and $h\in\Sch(s+\ima\R)$ for all $s\in\R$, then for each $s\in\R$
  there is $f_s\in\Sch(\R\inv_+)_s$ with $h(s+\ima t)=\widehat{f_s}(s+\ima t)$
  for all $t\in\R$.  Using the Cauchy--Riemann differential equation for the
  analytic function $h(s+\ima t)$, we conclude that~$f_s$ is independent
  of~$s$.  Hence we get a function~$f$ in~$\Hilm_-$ with $\hat{f}=h$ on all
  of~$\C$.  This yields the desired description of $(\Hilm_-)\sphat$.
  
  The same argument shows that $\Sch(\R\inv_+)_I\sphat$ for an open
  interval~$I$ is the space of all holomorphic functions $h\colon
  I+\ima\R\to\C$ with $h\in\Sch(s+\ima\R)$ for all $s\in I$.
  
  Proposition~\ref{pro:Sch_estimates} asserts that $f\in\Hilm_+$ if and only
  if $f\in\Sch(\R\inv_+)_{\ooival{0,\infty}}$ and $J\Fourier
  f\in\Sch(\R\inv_+)_{\ooival{-\infty,1}}$.  Since $J\Fourier$ is invertible
  on $\Sch(\R\inv_+)_s$ for $0<s<1$, this is equivalent to
  $f\in\Sch(\R\inv_+)_s$ for $s\ge½$ and $J\Fourier f\in\Sch(\R\inv_+)_s$
  for $s\le½$.  Moreover, these two conditions are equivalent for $s=½$.
  We get $f\in Z\Hilm_+$ if and only if $Z^{-1}f\in\Sch(\R\inv_+)_s$ for $s\ge
  ½$ and $J\Fourier Z^{-1}f\in\Sch(\R\inv_+)_s$ for $s\le ½$.
  
  Now let $h\colon\C\to\C$ be an entire function with $h\in\Sch(s+\ima\R)$ for
  all $s\in\R$.  Thus $h=\hat{f}$ for some $f\in\Hilm_-$.
  Equation~\eqref{eq:Zspectral} implies $f\in Z\Hilm_\cap$ if and only if
  $h/\zeta\in\Sch(s+\ima\R)$ for $s\ge½$ and $(J\Fourier f)\sphat/\zeta\in
  \Sch(s+\ima\R)$ for $s\le½$.  The functional
  equation~\eqref{eq:functional_equation} yields $(J\Fourier
  f)\sphat(z)/\zeta(z)=h(z)/\zeta(1-z)$.
\end{proof}

Let $(\Hilm^0_-)'$ be the space of continuous linear functionals
$\Hilm^0_-\to\C$ and let $\transpose{D}_-\in\End((\Hilm^0_-)')$ be the transpose
of~$D_-$.

\begin{corollary}  \label{cor:Dminus_eigenvalues}
  The eigenvalues of the transpose $\transpose{D}_-\in\End((\Hilm^0_-)')$ are
  exactly the zeros of the complete Riemann
  $\zeta$\nb-function~$\xi$.  The algebraic multiplicity
  of~$\lambda$ as an eigenvalue of $\transpose{D}_-$ is the zero order
  of~$\xi$ at~$\lambda$.  Here we define the algebraic multiplicity as
  the dimension of $\bigcup_{k\in\N} \ker (\transpose{D}_--\lambda)^k$.
\end{corollary}

\begin{proof}
  We may identify $(\Hilm^0_-)'$ with the space of continuous
  linear functionals on~$\Hilm_-$ that annihilate
  $Z\Hilm_\cap$.  Let $l\colon \Hilm_-\to\C$ be such a linear
  functional.  In the Fourier--Laplace transformed picture, we
  have $D_-h(s)=s\cdot h(s)$ for all $s\in\C$.  Hence~$l$ is an
  eigenvector for the transpose $\transpose{D}_-$ if and only
  if it is of the form $h\mapsto h(s)$ for some $s\in\C$.
  Similarly, $(\transpose{D}_--s)^k(l)=0$ if and only if
  $l(h)=\sum_{j=0}^{k-1} a_jh^{(j)}(s)$ for some
  $a_0,\dotsc,a_{k-1}\in\C$.  It follows from the functional
  equation $\xi(s)=\xi(1-s)$ that the zero orders of~$\xi$ at
  $s$ and $1-s$ agree for all $s\in\C$.  Moreover, $\xi$
  and~$\zeta$ have the same zeros and the same zero orders for
  $\RE s>0$.  Hence the assertion follows from
  Theorem~\ref{the:Lap_range}.
\end{proof}

\section{A geometric character computation}
\label{sec:geometric_character}

Tangermann: Consider case \(f=f_0* Jf_0\)!

Our next goal is to prove that the representation~$\rho$ is summable and to
compute its character geometrically.

\begin{definition}[\cite{Meyer:Primes_Rep}]
  Let~$G$ be a Lie group and let $\CCINF(G)$ be the space of smooth, compactly
  supported functions on~$G$.  A smooth representation~$\rho$ of~$G$ on a
  Fréchet space is called \emph{summable} if $\IN\rho(f)$ is nuclear for all
  $f\in\CCINF(G)$ and if these operators are uniformly nuclear for~$f$ in a
  bounded subset of $\CCINF(G)$.
\end{definition}

The theory of nuclear operators due to Alexandre
Grothendieck~\cite{Grothendieck:Produits} is rather
deep.  Nuclear operators are analogues of trace class operators on Hilbert
spaces.  It follows easily from the definition that $f\circ g$ is nuclear if
at least one of the operators $f$ and~$g$ is nuclear.  That is, the nuclear
operators form an \emph{operator ideal}.  For the purposes of this article, we
do not have to recall the definition of nuclearity because of the following
simple criterion:

\begin{theorem}  \label{the:nuclear_criterion}
  An operator between nuclear Fréchet spaces is nuclear if and only if it may
  be factored through a Banach space.
\end{theorem}

The spaces $\Sch(\R\inv_+)_I$ and~$\Hilm_+$ are nuclear because $\Sch(\R)$ is
nuclear and nuclearity is hereditary for subspaces and inverse limits.  Hence
Theorem~\ref{the:nuclear_criterion} applies to all operators between these
spaces.

The \emph{character of a summable representation}~$\rho$ is the distribution
on~$G$ defined by $\chi_\rho(f)\defeq \tr \IN\rho(f)$ for all $f\in\CCINF(G)$.
The uniform nuclearity of $\IN\rho(f)$ for~$f$ in bounded subsets of
$\CCINF(G)$ ensures that~$\chi_\rho$ is a bounded linear functional on
$\CCINF(G)$.  This is equivalent to continuity because $\CCINF(G)$ is an
LF-space.  If~$\rho$ is a virtual representation as in our case, we let $\tr
\IN\rho(f)$ be the supertrace $\tr\IN\rho_+(f)-\tr\IN\rho_-(f)$.

The above arguments and definitions show that summability of
representations really has to do with bounded subsets of $\CCINF(G)$ and
bounded maps, not with open subsets and continuous maps.  The same is true for
the concept of a nuclear operator.  That is, the theory of nuclear operators
and summable representations is at home in bornological vector spaces.  We may
still give definitions in the context of topological vector spaces if we turn
them into bornological vector spaces using the standard bornology of (von
Neumann) bounded subsets.  Nevertheless, topological vector spaces are the
wrong setup for studying nuclearity.  The only reason why I use them here is
because they are more familiar to most readers and easier to find in the
literature.

We need uniform nuclearity because we want $\chi_\rho(f)$ to be a bounded
linear functional of~$f$.  In the following, we will only prove nuclearity of
various operators.  The same proofs yield uniform nuclearity as well.  We
leave it to the reader to add the remaining details.  Suffice it to say that
there are analogues of Theorem~\ref{the:nuclear_criterion} and
Theorem~\ref{the:nuclear_hereditary} below for uniformly nuclear sets of
operators.

In order to prove the summability of our spectral interpretation~$\rho$, we
define several operators between~$\Hilm_\pm$ and the space
$\Sch_>\oplus\Sch_<$.  As auxiliary data, we use a smooth function
$\phi\colon\R_+\to[0,1]$ with $\phi(t)=0$ for $t\ll1$ and $\phi(t)=1$ for
$t\gg1$.  Let~$M_\phi$ be the operator of multiplication by~$\phi$.
We assume for simplicity that $\phi(t)+\phi(t^{-1})=1$, so that
\[
M_\phi+JM_\phi J=\ID.
\]
It is easy to check that~$M_\phi$ is a continuous map from~$\Sch_>$ into
$\Hol(\R\inv_+)$.  We warn the reader that our notation differs from that
in~\cite{Meyer:Primes_Rep}: there the auxiliary function $1-\phi$ is used and
denoted~$\phi$.
\begin{figure}[htb]
  \centering
  \setlength{\unitlength}{1cm}
  \begin{picture}(8,3)(0.2,0.2)
    \put(0.2,0.2){\vector(1,0){7.5}}
    \put(0.2,0.2){\vector(0,1){2.5}}
    \put(0.15,2.2){\line(1,0){0.1}}
    \put(2.7,0.15){\line(0,1){0.1}}
    \put(-0.1,2.07){$1$}
    \put(2.613,-0.2){$1$}
    \thicklines
    \put(4.2,2.2){\line(1,0){3.4}}
    \put(0.2,0.2){\line(1,0){1}}
    \qbezier(1.2,0.2)(2.2,0.2)(2.7,1.2)
    \qbezier(2.7,1.2)(3.2,2.2)(4.2,2.2)
  \end{picture}
  \caption{The auxiliary function~$\phi$}
  \label{fig:phi_picture}
\end{figure}
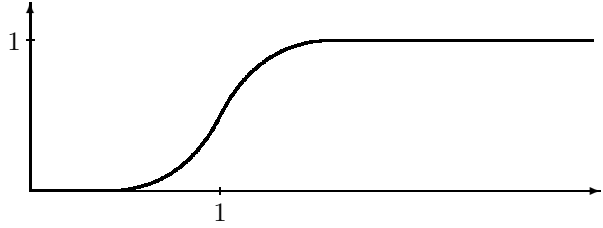

Now we define our operators:
\begin{alignat*}{2}
  \iota_+&\colon \Hilm_+\to\Sch_>\oplus\Sch_<,  &\qquad
  \iota_+ f &\defeq (Zf,JZ\Fourier f);
  \\
  \iota_-&\colon \Hilm_-\to\Sch_>\oplus\Sch_<,  &\qquad
  \iota_- f &\defeq (f,f);
  \\
  \pi_+&\colon \Sch_>\oplus\Sch_<\to\Hilm_+,  &\qquad
  \pi_+(f_1,f_2) &\defeq (M_\phi Z^{-1}f_1,\Fourier M_\phi Z^{-1} J f_2);
  \\
  \pi_-&\colon \Sch_>\oplus\Sch_<\to\Hilm_-,  &\qquad
  \pi_-(f_1,f_2) &\defeq (M_\phi f_1,J M_\phi J f_2).
\end{alignat*}
It follows from Equation~\eqref{eq:J_reflects} and
Proposition~\ref{pro:Z_invertible} that these operators are well-defined and
continuous.  The operators~$\iota_\pm$ are $\lambda$\nb-equivariant, the
operators~$\pi_\pm$ are not.  We compute
\begin{align*}
  \pi_-\iota_- &= M_\phi + JM_\phi J = \ID_{\Hilm_-},
  \\
  \pi_+\iota_+ &= M_\phi + \Fourier M_\phi\Fourier
  = \ID_{\Hilm_+} + M_\phi - \Fourier JM_\phi J\Fourier,
  \\
  \iota_-\pi_- &=
  \begin{pmatrix}
    M_\phi & JM_\phi J \\ M_\phi & JM_\phi J
  \end{pmatrix},
  \\
  \iota_+\pi_+ &=
  \begin{pmatrix}
    Z M_\phi Z^{-1} & Z\Fourier M_\phi Z^{-1}J \\
    JZ\Fourier M_\phi Z^{-1} & JZ M_\phi Z^{-1} J
  \end{pmatrix}.
\end{align*}
Thus~$\pi_-$ is a section for~$\iota_-$ and
$\iota_-\pi_-\in\End(\Sch_>\oplus\Sch_<)$ is a projection onto a subspace
isomorphic to~$\Hilm_-$.  The proof of Theorem~\ref{the:Zeta_estimate} shows
that~$\iota_+$ has closed range and is a topological isomorphism onto its
range.  Although~$\pi_+$ is not a section for~$\iota_+$, it is a near enough
miss for the following summability arguments.  I do not know whether there
is an honest section for~$\iota_+$, that is, whether the range
of~$\iota_+$ is a complemented subspace of $\Sch_>\oplus\Sch_<$.

\begin{lemma}  \label{lem:basic_nuclearity}
  The operator $\IN\lambda(f)(M_\phi - ZM_\phi Z^{-1})$ on~$\Sch_>$ is nuclear
  for $f\in\Sch_>$.

  The operator $\IN\lambda(f) M_\phi\colon \Sch_>\to\Sch_<$ is nuclear for
  $f\in\Sch_<$.

  The operator $\IN\lambda(f)( M_\phi - \Fourier JM_\phi J\Fourier)$
  on~$\Hilm_+$ is nuclear for $f\in\Hol(\R\inv_+)$.
\end{lemma}

\begin{proof}
  We have $\IN\lambda(f) Z(h) = f\ast \check{\zeta}\ast h = (Zf)\ast h =
  \IN\lambda(Zf)(h)$ for all $f,h\in\Sch_>$.  Hence
  \begin{multline*}
    \IN\lambda(f)(M_\phi - ZM_\phi Z^{-1})
    = [\IN\lambda(f),M_\phi] - [\IN\lambda(f) Z,M_\phi] Z^{-1}
    \\ = [\IN\lambda(f),M_\phi] - [\IN\lambda(Zf),M_\phi] Z^{-1},
  \end{multline*}
  We are going to show that $[\IN\lambda(f),M_\phi]$ is a nuclear operator
  on~$\Sch_>$ for all $f\in\Sch_>$.  Together with the above computation, this
  implies the first assertion of the lemma because~$Z^{-1}$ is continuous
  on~$\Sch_>$ and $f\in\Sch_>$ implies $Zf\in\Sch_>$.

  We compute
  \[
  [\IN\lambda(f),M_\phi] (h)(x) =
  \int_0^\infty f(xy^{-1}) \bigl(\phi(y)-\phi(x)\bigr) h(y) \,\diff\inv y.
  \]
  Thus our operator has the smooth integral kernel
  $f(xy^{-1})\bigl(\phi(y)-\phi(x)\bigr)$.  If~$f$ had compact support, this
  integral kernel would also be compactly supported.  If only $f\in\Sch_>$, we
  may still estimate that our kernel lies in $\Sch_>\hot \Sch(\R\inv_+)_{-s}$
  for any $s\in\ooival{1,\infty}$.  Thus $[\IN\lambda(f),M_\phi]$ factors
  through the embedding $\Sch_>\subseteq L^2(\R\inv_+,x^{2s}\,\diff\inv x)$ for
  any $s>1$.  This implies nuclearity by Theorem~\ref{the:nuclear_criterion}
  and finishes the proof of the first assertion of the lemma.
  
  To prove the second assertion, we let $L^2_I = \bigcap_{s\in
    I} L^2(\R\inv_+,x^{2s}\,\diff\inv x)$.  This is a Fréchet
  space.  We claim that $\IN\lambda(f)$ is a continuous linear
  operator $L^2_I\to\Sch_I$ for any open interval~$I$.
  (Actually, $\Sch_I$ is the Gårding subspace or, equivalently,
  the subspace of smooth vectors in the
  representation~$\lambda$ on~$L^2_I$, see
  also~\cite{Meyer:Smooth}.)  This follows from the description
  of~$\Sch_I$ in the proof of
  Proposition~\ref{pro:Sch_estimates} and $D^m(f*h)=(D^mf)*h$.
  Therefore, we have continuous linear operators
  \[
  \Sch_>\boin L^2_s \overset{M_\phi}\longrightarrow L^2_{\ocival{-\infty,s}}
  \boin L^2_{\ooival{-\infty,0}}
  \overset{\IN\lambda(f)}\longrightarrow \Sch_<
  \]
  for any $s>1$.  Thus $\IN\lambda(f)M_\phi$ factors through the Hilbert
  space~$L^2_s$.  This yields the assertion by our criterion for nuclear
  operators, Theorem~\ref{the:nuclear_criterion}.
  
  We claim that $\IN\lambda(f) (M_\phi - \Fourier J M_\phi J\Fourier)$ as an
  operator on~$\Hilm_+$ factors continuously through $L^2_{½} \cong
  L^2(\R\inv_+,x\,\diff\inv x)$.  This implies the third assertion.  To prove the
  claim, we use $\Hilm_+ = \Sch_>\cap \Fourier J(\Sch_<)$.  That is, the map
  \[
  \Hilm_+\to \Sch_>\oplus\Sch_<,
  \qquad f\mapsto (f,J\Fourier f)
  \]
  is a topological isomorphism onto its range.  We have already seen this
  during the proof of Theorem~\ref{the:Zeta_estimate}.  Hence we merely have
  to check the existence of continuous extensions
  \[
  \IN\lambda(f) (M_\phi - \Fourier J M_\phi J\Fourier)
  \colon L^2_{½}\to \Sch_>,
  \quad
  J\Fourier\circ \IN\lambda(f) (M_\phi - \Fourier J M_\phi J\Fourier)
  \colon L^2_{½}\to \Sch_<.
  \]
  We write
  \begin{multline*}
    \IN\lambda(f) (M_\phi - \Fourier J M_\phi J\Fourier)
    = -\IN\lambda(f) (M_{1-\phi} - \Fourier J M_{1-\phi} J\Fourier)
    \\ = -\IN\lambda(f) M_{1-\phi}
    + \IN\lambda(\Fourier Jf) M_{1-\phi} J\Fourier.
  \end{multline*}
  Proposition~\ref{pro:Sch_estimates} yields $\Fourier Jf\in\Sch_>$.
  Moreover, $J$ and~$\Fourier$ are unitary on $L^2_{½}$.  Hence we get
  bounded extensions $L^2_{½}\to \Sch_>$ of both summands by the same
  argument as for the second assertion of the lemma.  Similarly, both summands
  in
  \[
  J\Fourier\circ\IN\lambda(f) (M_\phi - \Fourier J M_\phi J\Fourier)
  = \IN\lambda(J\Fourier f) M_\phi - \IN\lambda(f) M_\phi J\Fourier
  \]
  have bounded extensions $L^2_{½}\to \Sch_<$ as desired.
\end{proof}

\begin{corollary}  \label{cor:almost_sections}
  The operator $\IN\lambda(f)\circ (\iota_-\pi_- - \iota_+\pi_+)$ on
  $\Sch_>\oplus\Sch_<$ is nuclear for all $f\in\Hol(\R\inv_+)$.
  
  The operator $\IN\lambda(f)\circ (\ID-\pi_+\iota_+)$ on~$\Hilm_+$ is nuclear
  for all $f\in\Hol(\R\inv_+)$.
\end{corollary}

\begin{proof}
  We have
  \begin{multline*}
    \IN\lambda(f)\circ (\iota_-\pi_- - \iota_+\pi_+)
    \\ =
    \begin{pmatrix}
      \IN\lambda(f)(M_\phi - Z M_\phi Z^{-1}) &
      \IN\lambda(f) J( M_\phi - JZ\Fourier M_\phi Z^{-1})J \\ 
      \IN\lambda(f)(M_\phi - JZ\Fourier M_\phi Z^{-1}) &
      \IN\lambda(f) J(M_\phi - Z M_\phi Z^{-1}) J
    \end{pmatrix}.
  \end{multline*}
  The upper left corner is nuclear by the first assertion of
  Lemma~\ref{lem:basic_nuclearity}.  Since
  \[
  \IN\lambda(f)(J(M_\phi - Z M_\phi Z^{-1}) J)
  = J(\IN\lambda(Jf)(M_\phi - Z M_\phi Z^{-1}))J,
  \]
  we also get the nuclearity of the lower right corner.  We have
  $\IN\lambda(f)JZ\Fourier= \IN\lambda(JZ\Fourier f)$ because $JZ\Fourier$ is
  $\lambda$\nb-invariant.  Proposition~\ref{pro:Sch_estimates} and
  Proposition~\ref{pro:Z_invertible} yield $JZ\Fourier(f)\in\Sch_<$.  Hence
  the two summands $\IN\lambda(f)M_\phi$ and $\IN\lambda(f)JZ\Fourier M_\phi
  Z^{-1}$ in the lower left corner are nuclear by the second assertion of
  Lemma~\ref{lem:basic_nuclearity}.  The assertion for the upper right corner
  follows by a symmetric argument.  The nuclearity of $\IN\lambda(f)\circ
  (\ID-\pi_+\iota_+)$ on~$\Hilm_+$ is exactly the third assertion of
  Lemma~\ref{lem:basic_nuclearity}.
\end{proof}

In order to apply this to the representation~$\rho$, we need a general fact
about nuclear operators.  Let $W_1,W_2$ be Fréchet spaces and let
$V_1\subseteq W_1$ and $V_2\subseteq W_2$ be closed subspaces.  If $T\colon
W_1\to W_2$ maps~$V_1$ into~$V_2$, we write $T|_{V_1,V_2}$ and
$T|^{W_1/V_1,W_2/V_2}$ for the operators $V_1\to V_2$ and $W_1/V_1\to W_2/V_2$
induced by~$T$.

\begin{theorem}[\cite{Grothendieck:Produits}]
  \label{the:nuclear_hereditary}
  If~$T$ is nuclear, so are $T|_{V_1,V_2}$ and $T|^{W_1/V_1,W_2/V_2}$.  If
  $V_1=V_2=V$ and $W_1=W_2=W$, then $\tr T= \tr T|_{V} + \tr T|^{W/V}$.
\end{theorem}

\begin{proposition}  \label{pro:rho_summable}
  The representation $\rho\colon \R\inv_+\to\Aut(\Hilm^0)$ is summable and
  \[
  \chi(\rho)(f) = \tr \IN\lambda(f) (\ID_{\Hilm_+} -\pi_+\iota_+)
  - \tr \IN\lambda(f) (\iota_-\pi_- -\iota_+\pi_+)
  \qquad
  \text{for all }f\in\CCINF(\R\inv_+).
  \]
\end{proposition}

\begin{proof}
  The embeddings of $\Hilm_+$ and~$\Hilm_-$ in $\Sch_>\oplus\Sch_<$ agree on
  the common subspace $\Hilm_\cap\cong Z\Hilm_\cap$ and hence combine to an
  embedding of~$\Hilm_\cup$.  Thus we identify~$\Hilm_\cup$ with the subspace
  $\iota_+\Hilm_+ + \iota_-\Hilm_-$ of $\Sch_>\oplus\Sch_<$.  Let
  $T\defeq\IN\lambda(f)(\iota_-\pi_- -\iota_+\pi_+)$.  The range of~$T$ is
  contained in~$\Hilm_\cup$.  Therefore, its trace as an operator on
  $\Sch_>\oplus\Sch_<$ agrees with its trace as an operator on~$\Hilm_\cup$ by
  Theorem~\ref{the:nuclear_hereditary}.  Write
  \[
  T|_{\Hilm_\cup}
  =
  \IN\lambda(f) (\ID_{\Hilm_\cup} -\iota_+\pi_+|_{\Hilm_\cup})
  - \IN\lambda(f) (\ID_{\Hilm_\cup} -\iota_-\pi_-|_{\Hilm_\cup}).
  \]
  Since $\iota_-\pi_-$ is a projection onto $\Hilm_-\subset\Hilm_\cup$,
  Theorem~\ref{the:nuclear_hereditary} yields
  \[
  \tr \IN\lambda(f) (\ID_{\Hilm_\cup} -\iota_-\pi_-|_{\Hilm_\cup})
  = \tr \IN\lambda(f)|^{\Hilm_\cup/\Hilm_-}
  = \chi(\rho_+)(f).
  \]
  Similarly, since $\iota_+\pi_+$ maps~$\Hilm_\cup$ into~$\Hilm_+$ we get
  \begin{multline*}
    \tr \IN\lambda(f) (\ID_{\Hilm_\cup} -\iota_+\pi_+|_{\Hilm_\cup})
    =
    \tr \IN\lambda(f)|^{\Hilm_\cup/\Hilm_+} +
    \tr \IN\lambda(f) (\ID_{\Hilm_+} -\iota_+\pi_+|_{\Hilm_+})
    \\ = \chi(\rho_-)(f) +
    \tr \IN\lambda(f) (\ID_{\Hilm_+} -\pi_+\iota_+).
  \end{multline*}
  Hence
  \[
  \tr \IN\lambda(f) (\iota_-\pi_- -\iota_+\pi_+)
  =
  \chi(\rho_-)(f) + \tr \IN\lambda(f) (\ID_{\Hilm_+} -\pi_+\iota_+) -
  \chi(\rho_+)(f).
  \]
  Along the way, we see that the operators whose trace we take are nuclear.
  That is, $\rho_+$ and~$\rho_-$ are summable representations.
\end{proof}

It remains to compute the traces in Proposition~\ref{pro:rho_summable}
explicitly.  We need the following definitions.  For a continuous function
$f\colon \R\inv_+\to\C$, let $\tau(f)\defeq f(1)$ and $\partial f(x)=f(x)\ln
x$.  This defines a bounded derivation~$\partial$ on $\Sch(\R\inv_+)_I$ for
any interval~$I$, that is, $\partial(f_1*f_2) =
\partial(f_1)*f_2+f_1*\partial(f_2)$.  This derivation is the generator of the
dual action $t\cdot f(x)\defeq x^{\ima t} f$ of~$\R$.  Notice that
$\tau(\partial f)=0$.  The obvious extension of~$\partial$ to distributions is
still a derivation.

\begin{lemma}  \label{lem:Toeplitz_trace}
  Let $f_0,f_1\in\Sch(\R\inv_+)_s$ for some $s\in\R$.  Then
  $\IN\lambda(f_0)[M_\phi,\IN\lambda(f_1)]$ is a nuclear operator on
  $L^2(\R\inv_+,x^{2s}\,\diff\inv x)$ and $\Sch(\R\inv_+)_s$ and
  \[
  \tr \IN\lambda(f_0)[M_\phi,\IN\lambda(f_1)]
  = \tau(f_0 * \partial f_1).
  \]
\end{lemma}

\begin{proof}
  The operators of multiplication by $x^{\pm s}$ are unitary operators between
  $L^2_s$ and $L^2_0$.  We may use them to reduce the general case to the
  special case $s=0$.  We assume this in the following.  We have checked above
  that $[M_\phi,\IN\lambda(f_1)]$ has an integral kernel in
  $\Sch(\R\inv_+)\hot\Sch(\R\inv_+)$.  Therefore, so has
  $\IN\lambda(f_0)[M_\phi,\IN\lambda(f_1)]$.  This implies nuclearity as an
  operator from $L^2(\R\inv_+,\diff\inv x)$ to $\Sch(\R\inv_+)$ by
  Theorem~\ref{the:nuclear_criterion}.  Moreover, the operator has the same
  trace on both spaces.  Explicitly, the integral kernel is
  \[
  (x,y)\mapsto \int_0^\infty f_0(xz^{-1}) f_1(zy^{-1})
  \bigl(\phi(z)-\phi(y)\bigr) \,\diff\inv z.
  \]
  We get
  \begin{multline*}
    \tr \IN\lambda(f_0)[M_\phi,\IN\lambda(f_1)]
    = \int_0^\infty \int_0^\infty f_0(xz^{-1}) f_1(zx^{-1})
    \bigl(\phi(z)-\phi(x)\bigr) \,\diff\inv z\,\diff\inv x
    \\ = \int_0^\infty f_0(x) f_1(x^{-1})
    \int_0^\infty \phi(z)-\phi(xz) \,\diff\inv z\,\diff\inv x.
  \end{multline*}
  We compute $\int_0^\infty \phi(z)-\phi(xz)\,\diff\inv z$.  If~$\phi$ had compact
  support, the $\lambda$\nb-invariance of $\diff\inv z$ would force the integral
  to vanish.  Therefore, we may replace~$\phi$ by any function~$\phi'$ with
  the same behaviour at $0$ and~$\infty$.  We choose~$\phi'$ to be the
  characteristic function of $\coival{1,\infty}$.  If $x\le1$, then
  $\phi'(z)-\phi'(xz)$ is the characteristic function of the interval
  $\coival{1,x^{-1}}$, so that the integral is $\ln(x^{-1})$.  We get the same
  value for $x\ge1$ as well.  Hence $\tr \IN\lambda(f_0)
  [M_\phi,\IN\lambda(f_1)] = \int_0^\infty f_0(x) f_1(x^{-1})
  \ln(x^{-1})\,\diff\inv x = \tau(f_0*\partial f_1)$.
\end{proof}

\begin{theorem}
  Define the distributions $W_p$ for $p\in\Primes$ and $p=\infty$ as in
  \eqref{eq:Wp} and~\eqref{eq:Winfty}.  Then
  \[
  \sum_{z\in\C} \ord_\xi(z) \hat{f}(z)
  = \chi(\rho)(f) = \sum_{p\in\Primes} W_p(f) + W_\infty(f)
  \]
  for all $f\in\Hol(\R\inv_+)$.  Here $\ord_\xi(z)$ denotes the order
  at~$z$ of the complete $\zeta$\nb-function~$\xi$, which is positive at
  poles and negative at zeros of~$\xi$.
\end{theorem}

\begin{proof}
  The trace of a nuclear operator on a nuclear Fréchet space is equal to the
  sum of its eigenvalues counted with algebraic multiplicity
  (see~\cite{Grothendieck:Produits}).  Since
  $\tr(A)=\tr({}^{t\!\!}A)$ for any nuclear operator~$A$, the first equality
  follows from Corollary~\ref{cor:Dminus_eigenvalues}.  It remains to show
  $\chi(\rho)(f)=\sum_{p\in\Primes} W_p(f) + W_\infty(f)$.
  Proposition~\ref{pro:rho_summable} yields
  \begin{multline}  \label{eq:get_W}
    \chi(\rho)(f)
    = - \tr \IN\lambda(f) (\iota_-\pi_- -\iota_+\pi_+)
    + \tr \IN\lambda(f) (\ID_{\Hilm_+} -\pi_+\iota_+)
    \\ = -\tr \IN\lambda(f) (M_\phi-ZM_\phi Z^{-1})|_{\Sch_>}
    - \tr \IN\lambda(f) J(M_\phi-ZM_\phi Z^{-1})J|_{\Sch_<}
    \\ {}- \tr \IN\lambda(f) (M_\phi-\Fourier JM_\phi J\Fourier)|_{\Hilm_+}.
  \end{multline}
  It suffices to check that this agrees with $W(f)$ if $f=f_0*f_1$ for
  $f_0,f_1\in\CCINF(\R\inv_+)$ because
  such elements are dense in $\Hol(\R\inv_+)$.  We compute
  \begin{multline*}
    - \tr \IN\lambda(f_0*f_1) (M_\phi-ZM_\phi Z^{-1})|_{\Sch_>}
    \\ = \tr \IN\lambda(f_0) [M_\phi,\IN\lambda(f_1)]
    - \tr \IN\lambda(f_0) [M_\phi,\IN\lambda(Zf_1)] Z^{-1}
    \\ = \tr \IN\lambda(f_0) [M_\phi,\IN\lambda(f_1)]
    - \tr \IN\lambda(Z^{-1} f_0) [M_\phi,\IN\lambda(Zf_1)]
  \end{multline*}
  because $\tr(AB)=\tr(BA)$ if~$A$ is nuclear.  This is a nuclear operator
  $L^2_s\to \Sch_>$ for any $s>1$.  Hence it has the same trace as an operator
  on $\Sch_>$ and $L^2_s$.  Lemma~\ref{lem:Toeplitz_trace} yields
  \begin{multline*}
    -\tr \IN\lambda(f_0*f_1) (M_\phi-ZM_\phi Z^{-1})|_{\Sch_>}
    = \tau(f_0 * \partial f_1) - \tau(Z^{-1}f_0 * \partial(Zf_1))
    \\ = - \tau(f_0 * f_1 * Z^{-1} \partial(Z))
    = \tau(f_0 * f_1 * Z \partial(Z^{-1}))
  \end{multline*}
  where $\partial(Z)$ is defined in the obvious way.  We also use
  $\partial(Z^{-1})=-Z^{-2} \partial(Z)$, which follows from the derivation
  property.  Now we use the Euler product for the Zeta operator and the
  derivation rule:
  \begin{multline*}
    Z*\partial(Z^{-1})
    = Z * \partial \biggl(\prod_{p\in\Primes} (1-\lambda^{-1}_p)\biggr)
    = \sum_{p\in\Primes} (1-\lambda^{-1}_p)^{-1} \partial(1-\lambda^{-1}_p)
    \\ = \sum_{p\in\Primes}  \ln(p) \lambda_p^{-1} (1-\lambda^{-1}_p)^{-1}
    =  \sum_{p\in\Primes} \sum_{e=1}^\infty \ln(p) \lambda_p^{-e}.
  \end{multline*}
  Hence
  \[
  -\tr \IN\lambda(f) (M_\phi-ZM_\phi Z^{-1})|_{\Sch_>}
  = \tau\biggl(f *\sum_{p\in\Primes} \sum_{e=1}^\infty \ln(p) \lambda_p^{-e}
  \biggr)
  = \sum_{p\in\Primes} \sum_{e=1}^\infty \ln(p) f(p^e).
  \]
  The second summand in~\eqref{eq:get_W} is reduced to this one by
  \begin{multline*}
    \tr \IN\lambda(f) J(M_\phi-ZM_\phi Z^{-1})J|_{\Sch_<}
    = \tr J\IN\lambda(f) J(M_\phi-ZM_\phi Z^{-1})|_{\Sch_>}
    \\ = \tr \IN\lambda(Jf) (M_\phi-ZM_\phi Z^{-1})|_{\Sch_>}.
  \end{multline*}
  Hence
  \[
  -\tr \IN\lambda(f) J(M_\phi-ZM_\phi Z^{-1})J|_{\Sch_<}
  =
  \sum_{p\in\Primes} \sum_{e=1}^\infty \ln(p) p^{-e} f(p^{-e}).
  \]
  These two summands together equal $\sum_{p\in\Primes} W_p$.
  
  Tangermann: The expression: $-\tau(f_0*f_1*J\Fourier \partial (FJ))$
  does not make sense exactly and seems a typo.  What does this mean
  if $f_0=f_1=\exp(-\pi x^2) p(x)$, where $p(x)$ is a polynomial
  satisfying certain vanishing conditions?

  Now we treat the third summand in~\eqref{eq:get_W}.  The same arguments as
  above yield\footnote{In the following computation, \(\Fourier J\) is viewed
    as a multiplier of \(\Sch(\R\inv_+)_{½}\) and~\(\partial\) is extended to
    multipliers using \(\partial(m)f\defeq \partial(mf)-m\partial(f)\).}
  \begin{multline*}
    -\tr \IN\lambda(f_0*f_1) (M_\phi-\Fourier JM_\phi J\Fourier)|_{\Hilm_+}
    \\ =  \tr \IN\lambda(f_0) [M_\phi,\IN\lambda(f_1)]
    - \tr \IN\lambda(f_0) [M_\phi,\IN\lambda(\Fourier Jf_1)] J\Fourier
    \\ = \tau(f_0*\partial f_1)
    - \tau(J\Fourier(f_0) * \partial (\Fourier Jf_1))
    = -\tau(f_0*f_1 * J\Fourier\partial (\Fourier J))
    = \tau(\Fourier J\partial (J\Fourier(f_0*f_1))).
  \end{multline*}
  Here we use $\tau(\partial f)=0$.  Explicitly,
  \[
  \tau(\Fourier J\partial (J\Fourier f))
  = \Fourier M_{\ln(x^{-1})} \Fourier^{-1} f(1)
  = -\Fourier(\ln x) \addconv f(1)
  = -\braket{\Fourier(\ln x)}{y\mapsto f(1-y)}.
  \]
  Here $M_{\ln(x^{-1})}$ denotes the operator of multiplication by
  $\ln(x^{-1})=-\ln x$ and~$\addconv$ denotes convolution with respect to the
  additive structure on~$\R$.  Thus it remains to compute the Fourier
  transform of $\ln x$.  Since $\psi\mapsto \int_\R \psi(x) \ln x\,\diff x$ defines
  a tempered distribution on~$\R$, $\Fourier(\ln x)$ is a well-defined
  tempered distribution on~$\R$.  The covariance property $\ln(tx)=
  \ln(t)+\ln(x)$ for $t,x\in\R\inv_+$ implies
  \begin{equation}  \label{eq:Fourier_ln_covariance}
    \braket{\Fourier(\ln x)}{\lambda_t\psi} =
    \braket{\Fourier(\ln x)}{\psi} - \ln(t)\psi(0).
  \end{equation}
  Especially, $\Fourier(\ln x)$ is $\lambda$\nb-invariant on the space of
  $\psi\in\Sch(\R)$ with $\psi(0)=0$.  Thus
  \[
  \braket{\Fourier(\ln x)}{\psi} = c\int_{\R\inv} \psi(x) \,\diff\inv x
  \]
  for some constant $c\in\R$ for all $\psi\in\Sch(\R)$ with $\psi(0)=0$.
  
  We claim that $c=-1$.  To see this, pick $\psi\in\Sch(\R)$ with
  $\psi(0)\neq0$ and consider $\psi-\lambda_t\psi$ for some $t\neq1$.
  Equation~\eqref{eq:Fourier_ln_covariance} yields
  \[
  c\int_{\R\inv} \psi(x) - \psi(t^{-1}x) \,\diff\inv x
  = \ln(t)\psi(0).
  \]
  As in the proof of Lemma~\ref{lem:Toeplitz_trace}, this implies $c=-1$.
  Thus the distribution $\Fourier(\ln x)$ is some principal value for the
  integral $-\int_{\R\inv} \psi(x) \abs{x}^{-1}\,\diff x$.  This principal value
  is described uniquely by the condition that $\Fourier^2(\ln x)(1)=0$.
  See also~\cite{Connes:Trace_Formula} for a comparison between this principal
  value and the one that usually occurs in the explicit formulas.
  
  Finally, we compute
  \begin{multline*}
    -\tr \IN\lambda(f) (M_\phi-\Fourier JM_\phi J\Fourier)|_{\Hilm_+}
    = -\braket{\Fourier(\ln x)}{y\mapsto f(1-y)}
    = \pv \int_{\R\inv} f(1-y) \frac{\diff y}{\abs{y}}
    \\ = \pv \int_{-\infty}^\infty f(x) \frac{\diff x}{\abs{1-x}}
    = \pv \int_0^\infty \frac{f(x)}{\abs{1-x}} + \frac{f(x)}{1+x} \,\diff x
    = W_\infty(f).
  \end{multline*}
  Plugging this into~\eqref{eq:get_W} gives the desired formula for
  $\chi(\rho)$.
\end{proof}

\section{Generalisation to Dirichlet \texorpdfstring{$L$}{L}-functions}
\label{sec:Dirichlet_L}

We recall the definition of Dirichlet $L$\nb-functions.  Fix some
$d\in\N_{\ge2}$ and let $(\Z/d\Z)\inv$ be the group of invertible elements in
the finite ring $\Z/d\Z$.  Let~$\chi$ be a character of $(\Z/d\Z)\inv$.
Define $\chi\colon \N\to\C$ by $\chi(n)\defeq\chi(n\bmod d)$ if $(n,d)=1$ and
$\chi(n)\defeq0$ otherwise.  The associated Dirichlet $L$\nb-function is
defined by
\[
L_\chi(s)
\defeq \sum_{n=1}^\infty \frac{\chi(n)}{n^s}.
\]
We suppose that~$d$ is equal to the conductor of~$\chi$, that is, $\chi$
does not factor through $(\Z/d'\Z)\inv$ for any proper divisor $d'\mid d$.  In
particular, $\chi\neq1$.

The constructions for the Riemann $\zeta$\nb-function that we have done
above work similarly for such $L$\nb-functions.  We define the
space~$\Hilm_-$ as above and let
\[
\Hilm_+ \defeq \{f\in\Sch(\R) \mid f(-x) = \chi(-1)f(x)\}
\]
be the space of even or odd functions, depending on $\chi(-1)\in\{\pm1\}$.
The Fourier transform on $\Sch(\R)$ preserves the subspace~$\Hilm_+$, and the
assertions of Proposition~\ref{pro:Sch_estimates} remain true.  However, now
$\Fourier^2=\chi(-1)$, so that we have to replace~$\Fourier$ by
$\Fourier^*=\chi(-1)\Fourier$ in appropriate places to get correct formulas.

Of course, the $L$\nb-function analogue of the Zeta operator is defined by
\[
\mathcal{L}_\chi f(x) \defeq \sum_{n=1}^\infty \chi(n)\cdot f(nx),
\]
for $f\in\Hilm_+$.  We now have the Euler product expansion
\[
\mathcal{L}_\chi
= \sum_{n=1}^\infty \chi(n) \lambda_{n^{-1}}
= \prod_{p\in\Primes} \sum_{e=0}^\infty \chi(p)^e \lambda_{p^{-e}}
= \prod_{p\in\Primes} (1-\chi(p)\lambda_{p^{-1}})^{-1}.
\]
The same estimates as for the Zeta operator show that this product
expansion converges on $\Sch_>$ (compare Proposition~\ref{pro:Z_invertible}).

The Poisson Summation Formula looks somewhat different now: we have
\[
\mathcal{L}_\chi(f)
= \kappa\cdot d^{½} \lambda_d^{-1} J\mathcal{L}_\chi\Fourier(f),
\]
where~$\kappa$ is some complex number with $\abs{\kappa}=1$.  The proof
of Theorem~\ref{the:Zeta_estimate} then carries over without change.  We also
get the holomorphic continuation and the functional equation for~$L_\chi$.  If
$\chi(-1)=-1$, we have to use the special function $2x\exp(-\pi x^2)$ instead
of $2\exp(-\pi x^2)$ to pass from~$L_\chi$ to the complete $L$\nb-function.
The results in Section \ref{sec:spectral_interpretation} carry over in an
evident way.  Now $\Hilm^0_+=\{0\}$ because~$L_\chi$ does not have poles, and
the eigenvalues of $\transpose{D}_-$ are the non-trivial zeros of~$L_\chi$, with
correct algebraic multiplicity.

Some modifications are necessary in Section~\ref{sec:geometric_character}.  We
define $\iota_-$ and~$\pi_-$ as above.  Since we want the
embeddings~$\iota_\pm$ to agree on $\Hilm_\cap\cong
\mathcal{L}_\chi\Hilm_\cap$, we should put
\[
\iota_+(f) \defeq (\mathcal{L}_\chi(f),
\kappa\cdot d^{½} \lambda_d^{-1} J\mathcal{L}_\chi\Fourier f)
\]
and modify~$\pi_+$ accordingly so that $\pi_+\iota_+=M_\phi +\Fourier^*
M_\phi\Fourier$.  With these changes, the remaining computations carry over
easily.  Of course, we get different local summands~$W_p$ in the explicit
formula for~$\mathcal{L}_\chi$.  You may want to compute them yourself as an
exercise to test your understanding of the arguments above.

Even more generally, we may replace the rational numbers~$\Q$ by an imaginary
quadratic extension like $\Q[\ima]$ and study $L$\nb-functions attached to
characters of the idèle class group of this field extension.  Such fields have
only one infinite place, which is complex.  A character of the idèle class
group restricts to a character of the circle group inside~$\C\inv$.  The
space~$\Hilm_+$ is now replaced by the homogeneous subspace of $\Sch(\C)$
defined by that character.

Once there is more than one infinite place, we need the more general setup
of~\cite{Meyer:Primes_Rep}.  In addition, the adèlic constructions
in~\cite{Meyer:Primes_Rep} provide a better understanding even for~$\Q$
because they show the similarity of the analysis at the finite and infinite
places.  The explicit formula takes a much nicer form if we put together all
characters of the idèle class group.  The resulting local summands that make
up the Weil distribution are of the same general form
\[
W_v(f) = \pv\int_{\Q_v\inv} \frac{f(x)}{\abs{1-x^{-1}}} \,\diff\inv x 
\]
at all places~$v$ (with $\Q_\infty=\R$) and may also be interpreted
geometrically as a generalised Lefschetz trace formula
(see~\cite{Connes:Trace_Formula}).

\end{document}